\documentclass[a4paper,11pt]{article}	
\usepackage[english]{babel}
\usepackage[utf8]{inputenc}

\usepackage{amsmath,amsthm,amssymb,graphicx}
\usepackage{hyperref,color,amsthm,mathdots}
\usepackage{enumerate}
\usepackage{cite}

\newtheorem{property}{Property}[section]

\newtheorem{proposition}[property]{Proposition}
\newtheorem{theorem}[property]{Theorem}
\newtheorem{lemma}[property]{Lemma}

\newtheorem{observation}[property]{Observation}
\newtheorem{claim}{Claim}

\theoremstyle{definition}

\newtheorem{remark}[property]{Remark}

\newcommand{\bw}{{\bf bw}}
\newcommand{\dist}{{\bf dist}}

\begin{document}

\title{Forbidding Kuratowski Graphs as Immersions}

\date{}
\author{Archontia C. Giannopoulou\thanks{Department of Mathematics, National and Kapodistrian University of 	Athens, Panepistimioupolis, GR-15784,  Athens, Greece. Emails: \texttt{\{arcgian,sedthilk\}@math.uoa.gr}}~\thanks{Supported by a grant of the Special Account for Research Grants of the  National and Kapodistrian University of Athens (project code: 70/4/10311).}
	\and
	Marcin Kamiński\thanks{Département d'Informatique, 
	Université Libre de Bruxelles, Brussels, Belgium. Email: \texttt{Marcin.Kaminski@ulb.ac.be}}
	\and
	Dimitrios M. Thilikos$^{*}$\thanks{Co-financed by the European Union (European Social Fund - ESF) and Greek national funds through the Operational Program ``Education and Lifelong Learning'' of the National Strategic Reference Framework (NSRF) - Research Funding Program: ``Thales. Investing in knowledge society through the European Social Fund.''}
}

\maketitle

\begin{abstract}\noindent 
The immersion relation is a partial ordering relation on graphs that is 
weaker than the topological minor relation in the sense that if a graph $G$ contains a 
graph $H$ as a topological minor, then it also contains it as an immersion 
but not vice versa. Kuratowski graphs, namely $K_{5}$ and $K_{3,3}$, 
give a precise characterization of planar graphs when excluded as 
topological minors.
In this note we give a structural characterization of the graphs that exclude Kuratowski graphs as immersions. 
 We prove that they can be constructed by applying consecutive $i$-edge-sums, for $i\leq 3$, starting from graphs that are planar sub-cubic or of branch-width at most 10.
\medskip

\end{abstract}
\section{Introduction}

A famous graph theoretic result is the theorem of Kuratowski which states that a graph $G$ is planar if and only if it does not contain $K_{5}$ and $K_{3,3}$ (also known as the Kuratowski graphs) as a topological minor, that is, if $K_{5}$ and $K_{3,3}$ cannot be obtained from the graph by applying vertex and edge removals and edge dissolutions.
It is well-known that the topological minor relation defines a (partial) ordering of the class of graphs. 

In a similar way, the immersion and the minor orderings can be defined in graphs if instead of vertex dissolutions we ask for edge lifts and edge contractions respectively. (For detailed definitions see Section~\ref{sec:defn}.) Notice here that the topological minor ordering is stronger than the minor and the immersion orderings in the sense that if a graph $G$ contains a graph $H$ as a topological minor then it also contains it as an immersion and as a minor but the inverse direction does not always hold.

In the celebrated Graph Minors theory, developed by Robertson and Seymour, it was proven that both the immersion and minor orderings are well-quasi-ordered, that is, there are no infinite sets of mutually non-comparable graphs~\cite{RobertsonS03,RobertsonS10} according to these orderings.
This result has as a consequence the complete characterization of the graph classes that are closed under taking immersions or minors in terms of forbidden graphs, where a graph class is closed under taking immersions (respectively minors) if for any graph that belongs to the graph class all of its immersions (respectively minors) also belong to the graph class.
For example, by an extension of the Kuratowski theorem (also known as Wagner's theorem), it is also known that a graph is planar if and only if it does not contain $K_{5}$ and $K_{3,3}$ as a minor.

Thus, a  question that naturally arises is about the characterization of the structure of a graph $G$ that excludes some fixed graph $H$ as an immersion or as a minor. While this subject has been extensively studied for the minor ordering (see~\cite{DvorakGT12,RueST12,RobertsonST94,RobertsonS-XIII,abs-1102-5762,ArnborgPC90,BodlaenderT99,KoutsonasTY11,Robertson93linklessembeddings,TakahashiUK94}), the immersion ordering only recently attracted the attention of the research community~\cite{GroheKMW11,2011arXiv1101.2630D,GiannopoulouSZ12,Abu-KhzamL03,KawarabayashiK12}.
In~\cite{2011arXiv1101.2630D}, DeVos et al. proved that for every positive integer $t$, every simple graph of minimum degree at least $200t$ contains the complete graph on $t$ vertices as a (strong) immersion
and in~\cite{FerraraGTW08} Ferrara et al., given a graph $H$, provide a lower bound on the minimum degree of any graph $G$ in order to ensure that $H$ is contained in $G$ as an immersion.
More recently, in~\cite{SeymourW2012}, Seymour and Wollan proved a structure theorem for graphs excluding complete graphs as immersions.

In terms of graph colorings, Abu-Khzam and Langston in~\cite{Abu-KhzamL03} provide evidence supporting the analog of Hadwiger's Conjecture according to the immersion ordering, that is, the conjecture stating that if the chromatic number of a graph $G$ is at least $t$ then $G$ contains the complete graph on $t$ vertices as an immersion and prove it for $t\leq 4$. This conjecture is proven for $t=5,6$ and $t\leq 7$ by  Lescure and Meyniel in~\cite{Lescure1988325} and by DeVos et al. in~\cite{1213.05137} independently.
The most recent result on colorings is an approximation of the list coloring number on graphs excluding the complete graph as immersion~\cite{KawarabayashiK12}.

Finally, in terms of algorithms,  in~\cite{GroheKMW11}, Grohe et al. gave a cubic time algorithm that decides whether a fixed graph $H$ is contained in an input graph $G$ as immersion and in~\cite{GiannopoulouSZ12} it was proved that the minimal graphs not belonging to a graph class closed under immersions can be computed when an upper bound on their tree-width and a description of the graph class in Monadic Second Order Logic are given.

In this note we characterize the structure of the graphs that do not contain $K_{5}$ and $K_{3,3}$ as immersions. 
As these graphs already exclude Kuratowski graphs as topological minors they are already planar.
Additionally, we show that they have a more special structure: they can be constructed by repetitively, 
joining together simpler graphs, starting from 
either graphs of small decomposability or by planar graphs with maximum degree 3.
In particular, we prove that a graph $G$ that does not contain neither $K_{5}$ nor $K_{3,3}$ as immersions can be constructed by applying consecutive $i$-edge-sums, for $i\leq 3$, to graphs that are planar sub-cubic or of branch-width at most 10.

\section{Definitions}\label{sec:defn}
For every integer $n$, we let $[n]=\{1,2,\dots,n\}$.
All graphs we consider are finite, undirected, and loopless but may have  multiple edges. 
Given a graph $G$ we denote by $V(G)$ and $E(G)$
its {\em vertex} and {\em edge set} respectively.
Given a set $F\subseteq E(G)$ (resp. $S\subseteq V(G)$), we denote 
by $G\setminus F$ (resp. $G\setminus S$) the graph obtained from $G$ if we remove the edges in $F$ (resp. the vertices in $S$ along with their incident edges). 
We denote by ${\cal C}(G)$ the set of the {\em connected components} of $G$.
Given  a vertex $v\in V(G)$, we also use the notation $G\setminus v=G\setminus \{v\}$. 
The {\em neighborhood} 
of a vertex $v\in V(G)$, denoted by $N_{G}(v)$, is the set of 
edges in $G$ that are adjacent to $v$. We denote by $E_{G}(v)$ the set of the edges 
of $G$ that are incident to $v$.
The {\em degree} of a vertex $v\in V(G)$, denoted by $\deg_{G}(v)$, 
is the number of edges that are incident to it, i.e., $\deg_{G}(v)=|E_{G}(v)|$. 
Notice that, as we are dealing with multigraphs,  
$|N_{G}(v)|\leq \deg_{G}(v)$. The minimum degree of a graph $G$, denoted by $\delta(G)$, is the minimum of the degrees of the vertices of $G$, that is, $\delta(G)=\min_{v\in V(G)}\deg_{G}(v)$.
A graph is called {\em sub-cubic} if all its vertices have degree at most 3. 
We  also denote by $K_{r}$ the {\em complete graph} on $r$ vertices and by $K_{r,q}$ the {\em complete 
bipartite graph} with $r$ vertices in its one part and $q$ in the other. 
Let $P$ be a path and $v,u\in V(P)$. We denote by $P[v,u]$ the sub-path of $P$ with end-vertices $v$ and $u$.
Given two paths $P_{1}$ and $P_{2}$ who share a common endpoint $v$, we say that they are {\em well-arranged} if their common vertices appear in the same order in both paths.

We say that a graph $H$ is a {\em subgraph} 
of a graph $G$ if $H$ can be obtained 
from $G$, after removing edges or vertices.  
An {\em edge cut} in a graph $G$ is a non-empty set $F$ of edges that belong to the same connected component 
of $G$ and such that $G\setminus F$ has 
more connected components than $G$. If $G\setminus F$ has one more 
connected component than $G$ then we say that $F$ is a {\em minimal} edge cut.
Let $F$ be an edge cut of a graph $G$ and let $G'$  be the connected component of $G$ containing 
the edges of  $F$. We say that $F$ is an  {\em internal edge cut} if it is minimal 
and both connected components of $G'\setminus F$
contain at least 2 vertices. An edge cut is also called {\em $i$-edge-cut} if it has cardinality  $\leq i$.

In this paper we mostly deal with lanai graphs, that is graphs that are embedded 
in the sphere $\Bbb{S}_{0}$. We call such a graph, along with its embedding, {\em $\Sigma_{0}$-embeddable 
graph}. Let $C_{1}$, $C_{2}$ be two disjoint cycles  in a  $\Sigma_{0}$-embeddable 
graph $G$. Let also $\Delta_{i}$ be the open disk of $\Bbb{S}_{0}\setminus C_{i}$ that does not contain points of $C_{3-i}$, $i\in [2]$. The {\em annulus between}  $C_{1}$ and $C_{2}$ is the  set $\Bbb{S}_{0}\setminus (\Delta_{1}\cup\Delta_{2})$ and we denote it by $A[C_{1},C_{2}]$. Notice that  $A[C_{1},C_{2}]$ is  a closed set.  
If ${\cal A}=\{C_{1},\ldots,C_{r}\}$ is a collection of cycles of a $\Bbb{S}_{0}$-embeddible graph $G$.
We say that ${\cal A}$ is {\em nested} if for every $i\in[r-2]$ $C_{i+1}$,  $A[C_{i},C_{i+1}]\cup A[C_{i+1}\cup A_{i+1}]=A[i,i+2]$.

\paragraph{Contractions and minors.}
The {\em contraction of 
an edge} $e=\{x,y\}$ from $G$ is the removal from 
$G$ of all edges incident to $x$ or $y$ and 
the insertion of a new vertex $v_{e}$ that is made adjacent 
to all the vertices of $(N_{G}(x)\setminus \{y\})\cup (N_{G}(y)\setminus \{x\})$ such 
that edges corresponding to the vertices in $(N_{G}(x)\setminus \{y\})\cap (N_{G}(y)\setminus \{x\})$
increase their multiplicity, that is, if there was a vertex $v\in (N_{G}(x)\setminus \{y\})\cap (N_{G}(y)\setminus \{x\})$, $k$ edges joining $v$ and $x$ and, $l$ edges joining $v$ and $y$ then in the resulting graph there will be $k+l$ edges joining $v$ with $v_{e}$. Finally, remove any loops resulting from this operation.
Given two graphs $H$ and $G$, we say that  $H$ is a {\em contraction}
of $G$, denoted by $H\leq_{c}G$, if $H$ can be obtained from $G$ after a (possibly empty) 
series of edge contractions. Moreover,
$H$ is a {\em minor} of $G$
if $H$ is a contraction of some subgraph of $G$.

\paragraph{Topological minors.} A {\em subdivision} of
a graph $H$ is any graph obtained after replacing some of its edges 
by paths between the same endpoints.
A graph $H$ is a {\em topological minor} of $G$ (denoted by $H\leq_{t} G$) 
if $G$ contains
as a subgraph some subdivision of $H$.

\paragraph{Immersions.}  The {\em lift} of two edges $e_{1}=\{x,y\}$
and $e_{2}=\{x,z\}$ to an edge $e$ is the operation of removing $e_{1}$ and $e_{2}$ from $G$ and then adding the edge $e=\{y,z\}$ in the resulting graph.
We say that a graph $H$ can be {\em (weakly) immersed} in a graph $G$ (or is an {\em immersion} of $G$), denoted by $H\leq_{im} G$,  if $H$ can 
be obtained from a subgraph of $G$ after a (possibly empty) sequence of edge lifts. Equivalently, we say that $H$ is an immersion of $G$ if there is an injective mapping $f: V(H) \to V(G)$ such that, for every edge $\{u,v\}$ of $H$, there is a path from $f(u)$ to $f(v)$ in $G$ and for any two distinct edges of $H$ the corresponding paths in $G$ are {\em edge-disjoint}, that is, they do not share common edges. 
Additionally, if these paths are internally disjoint from $f(V(H))$, then we say that $H$ is {\it strongly immersed} in $G$ (or is a {\em strong immersion of $G$}).
The injective mapping $f$ together with the edge-disjoint paths is called {\em a model of $H$ in $G$ defined by $f$}.

\begin{figure}[h]
\begin{center}
\scalebox{.75}{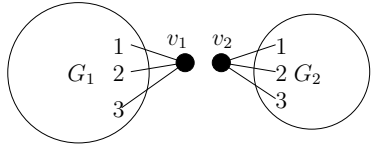} 
\end{center}
\caption{The graphs $G_{1}$ and $G_{2}$ before the edge-sum}
\label{fig:edgsmbef}
\end{figure}

\begin{figure}[h]
\begin{center}
\scalebox{.75}{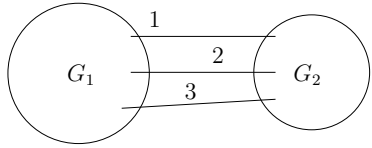} 
\end{center}
\caption{The graph obtained after the edge-sum}
\label{fig:edgsmaf}
\end{figure}

\paragraph{Edge sums.} 
Let $G_{1}$ and $G_{2}$ be graphs, let $v_{1},v_{2}$ be vertices of $V(G_{1})$ and $V(G_{2})$ respectively, 
and consider a bijection $\sigma: E_{G_{1}}(v_{1})\rightarrow E_{G_{2}}(v_{2})$ where $E_{G_{1}}(v_{1})=\{e_{1}^{i}\mid i\in [k]\}$.  We define the {\em $k$-edge sum} of 
$G_{1}$ and $G_{2}$ on $v_{1}$ and $v_{2}$ as the graph $G$ obtained if we take 
the disjoint union of $G_{1}$ and $G_{2}$, identify $v_{1}$ with $v_{2}$, and
then,  for each $i\in\{1,\ldots,k\}$, lift $e_{1}^{i}$ and $\sigma(e_{1}^{i})$ to 
a new edge $e^{i}$ and remove the vertex $v_{1}$. (See Figures~\ref{fig:edgsmbef} and~\ref{fig:edgsmaf})

Let $G$ be a graph, let $F$ be a minimal $i$-edge cut in $G$, and let $G'$ be the connected component 
of $G$ that contains $F$. Let also $C_{1}$ and $C_{2}$ be the two  connected components of 
$G'\setminus F$. We denote by $C'_{i}$ the graph obtained from $G'$ 
after contracting all edges of $C'_{3-i}$ to a single vertex $v_{i}, i\in\{1,2\}$.  
We say that the graph consisting of the disjoint union of the graphs in ${\cal C}(G)\setminus \{C_{1},C_{2}\}\cup\{C'_{1},C'_{2}\}$ is the  {\em $F$-split} 
of $G$ and we denote it by  $G\!\!\mid_{F}$. 
Notice that if $G$ is connected and $F$ is a minimal $i$-edge cut in $G$, then 
$G$ is the result of the $i$-edge sum of the two connected components $G_{1}$ and $G_{2}$
of ${\cal C}(G\!\!\mid_{F})$ on the vertices $v_{1}$ and $v_{2}$. From Menger's Theorem we obtain the following.

\begin{observation}\label{obs:mngr}
Let $k$ be a positive integer. If $G$ is a connected graph that does not contain an internal $i$-edge cut, for some $i\in [k-1]$ and $v,v_{1},\dots, v_{i}\in V(G)$ are distinct vertices such that $\deg_{G}(v)\geq i$ then there exist $i$ edge-disjoint paths from $v$ to $v_{1},v_{2},\dots,v_{i}$.
\end{observation}

\begin{lemma}
\label{frloo}
If $G$ is a $\{K_{5},K_{3,3}\}$-immersion free  
connected graph and $F$ is a minimal internal $i$-edge cut in $G$, for $i\in [3],$
then both connected components of $G\!\!\mid_{F}$ are  $\{K_{5},K_{3,3}\}$-immersion free.
\end{lemma}

\begin{proof}
For contradiction assume that $G$ is a $\{K_{5},K_{3,3}\}$-immersion free  
connected graph and one of the connected components of $G\!\!\mid_{F}$, say $C_{1}'$, contains $K_{5}$ or $K_{3,3}$ as an immersion, where $F$ is a minimal internal $i$-edge cut in $G$, $i\in [3]$. 
Assume that $H\in\{K_{5},K_{3,3}\}$ is immersed in $C_{1}'$ and let $f:V(H)\rightarrow V(C_{1}')$ be a model of $H$ in $G$. Let also $v_{1}$ be the newely introduced vertex of $C_{1}'$.
Notice that if $v_{1}\notin f(V(H))$ and $v_{1}$ is not an internal vertex of any of the edge-disjoint paths between the vertices in $f(V(H))$, then $f$ is a model of $H$ in $C_{1}$. As $C_{1}\subseteq G$, $f$ is a model of $H$ in $G$, a contradiction to the hypothesis. Thus, we may assume that either $v_{1}\in f(V(H))$ or $v_{1}$ is an internal vertex in at least one of the edge-disjoint paths between the vertives in $V(H)$.
Note that, as neither $K_{5}$ nor $K_{3,3}$ contain vertices of degree $1$, $\deg_{C_{1}'}(v_{1})= 2$ or $\deg_{C_{1}'}(v_{1})=3$. 

We first exclude the case where $v_{1}\notin f(V(H))$, that is, $v_{1}$ only appears as an internal vertex on the edge-disjoint paths. Observe that, as $\deg_{C_{1}'}(v_{1})\leq 3$, $v_{1}$ belongs to exactly one path $P$ in the model defined by $f$. Let $v_{1}^{1}$ and $v_{1}^{2}$ be the neighbors of $v_{1}$ in $P$. Recall that, by the definition of an internal $F$-split, there are vertices $v_{2}^{1}$ and $v_{2}^{2}$ in $C_{2}$ such that $\{v_{1}^{1},v_{2}^{1}\},\{v_{2}^{1},v_{2}^{2}\}\in E(G)$. Furthermore, as $C_{2}$ is connected, there exists a $(v_{2}^{1},v_{2}^{2})$-path $P'$ in $C_{2}$. Therefore, be substituting the subpath $P[v_{1}^{1},v_{1}^{2}]$ by the path defined by the union of the edges $\{v_{1}^{1},v_{2}^{1}\},\{v_{2}^{1},v_{2}^{2}\}\in E(G)$ and the path $P'$ in $C_{2}$ we obtain a model of $H$ in $G$ defined by $f$, a contradiction to the hypothesis.

Thus, the only possible case is that $v_{1}\in f(V(H))$. As $\delta(K_{5})=4$ and $\deg_{C_{1}'}(v_{1})\leq 3$, $f$ defines a model of $K_{3,3}$ in $C_{1}'$. Let $v_{1}^{1},v_{1}^{2}$ and $v_{1}^{3}$ be the neighbors of $v_{1}$ in $C_{1}'$. We claim that there is a vertex $v$ in $C_{2}$ and edge-disjoint paths from $v$ to $v_{1}^{1},v_{1}^{2},v_{1}^{3}$ in $G$, thus proving that there exists a model of $K_{3,3}$ in $G$ as well, a contradiction to the hypothesis. By the definition of an internal $F$-split, there are vertices $v_{2}^{1},v_{2}^{2}$ and $v_{2}^{3}$ in $C_{2}$ such that $\{v_{1}^{i},v_{2}^{i}\}\in E(G)$, $i\in [3]$.
Recall that $C_{2}$ is connected. Therefore, if for every vertex $v\in C_{2}$, $\deg_{C_{2}}(v)\leq 2$, $C_{2}$ contains a path whose endpoints, say $u$ and $u'$ belong to $\{v_{2}^{1},v_{2}^{2},v_{2}^{3}\}$ and internally contains the vertex in $\{v_{2}^{1},v_{2}^{2},v_{2}^{3}\}\setminus \{u,u'\}$, say $u''$. Then it is easy to verify that $u''$ satisfies the conditions of the claim. 
Assume then that there is a vertex $v\in C_{2}$ of degree at least 3. Let $G'$ be the graph obtained from $G$ after removing all vertices in $V(C_{1})\setminus \{v_{1}^{1},v_{1}^{2},v_{1}^{3}\}$ and adding a new vertex that we make it adjacent to the vertices in $\{v_{1}^{1},v_{1}^{2},v_{1}^{3}\}$. As $G$ does not contain an internal $i$-edge cut, $i\in [2]$, $G'$ does not contain an internal $i$-edge cut, $i\in [2]$. Therefore, from Observation~\ref{obs:mngr} and the fact that $v\notin \{v_{1}^{1},v_{1}^{2},v_{1}^{3}\}$, we obtain that there exist 3 edge-disjoint paths from $v$ to $v_{1}^{1},v_{1}^{2},v_{1}^{3}$ in $G'$ and thus in $G$. This completes the proof of the claim and the lemma follows.  
\end{proof}

\begin{figure}
\label{cyl44}
\begin{center}
\scalebox{0.35}{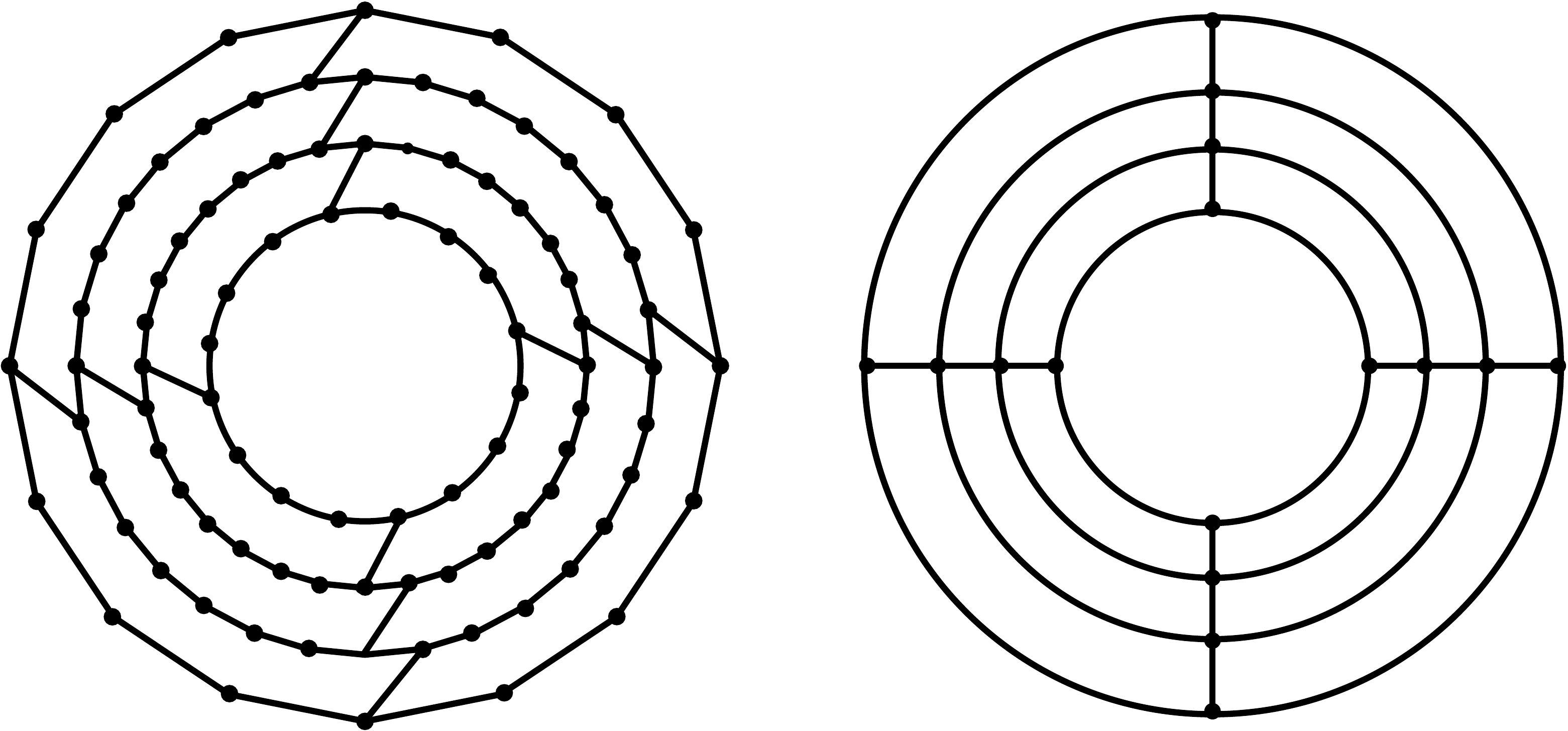}
\end{center}
\caption{A (4,4)-railed annulus and a (4,4)-cylinder}
\end{figure}

Let $r \geq 3$ and $q\geq 1$. A $(r,q)$-{\em cylinder}, denoted by $C_{r,q}$, is the cartesian product of a cycle on $r$ vertices and a path on $q$ vertices. (See, for example, Figure~\ref{cyl44}) 
A $(r,q)$-{\em railed annulus} in a graph $G$ is a pair $({\cal A},{\cal W})$ such that ${\cal A}$ is a collection of $r$ nested cycles $C_{1},C_{2},\ldots ,C_{r}$ that are all met by a collection ${\cal W}$ of $q$ paths $P_{1},P_{2},\dots, P_{q}$ (called {\em rails}) in a way that the intersection of a rail and a path is always a (possibly trivial, that is, consisting of only one vertex) path. (See, for example, Figure~\ref{cyl44})
Notice that given a graph $G$ embedded in the sphere and a $(k,h)$-{\em cylinder} ($(r,q)$-{\em railed annulus} respectively) of $G$, then any two cycles of the $(k,h)$-{\em cylinder} ($(r,q)$-{\em railed annulus} respectively) define an annulus between them.

\paragraph{Branch decompositions.}
A {\em branch decomposition of a graph} $G$ is a pair $B=(T,\tau)$, where $T$ is a ternary tree and $\tau: E(G)\rightarrow {\cal L}(T)$ is a bijection of the edges of $G$ to the leaves of ${T}$, denoted by ${\cal L}(T)$. Given a branch decomposition ${\cal B}$, we define $\sigma_{B}:E(T)\rightarrow \mathbb{N}$ as follows.

Given an edge $e\in E(T)$ let $T_{1}$ and $T_{2}$ be the trees in $T\setminus \{e\}.$ Then $\sigma_{B}(e)=|\{v\mid \text{there exist } e_{i}\in \tau^{-1}({\cal L}(T_{i})) \text{, } i\in[2] \text{, such that } e_{1}\cap e_{2}=\{v\}\}|$.
 The {\em width of a branch decomposition} $B$ is $\max_{e\in E(T)}\sigma_{B}(e)$ and the {\em branch-width of a graph} $G$, denoted by $\bw(G)$, is the minimum width over all branch decompositions of $G$. In the case where $|V(T)|\leq 1$ the width of the branch decomposition is defined to be 0.
The following has been proven in~\cite{GuT10}.

\begin{theorem}\label{thm:bnwdthcyl}
If $G$ is a planar graph and $k,h$ are integers with $k\geq 3$ and $h\geq 1$ then $G$ either contains the $(k,h)$-cylinder as a minor or has branch-width at most $k+2h-2$. 
\end{theorem}

We now prove the following.

\begin{lemma}\label{lem:brnwdthranl}
If $G$ is a planar graph of branch-width at least 11, then $G$ contains a (4,4)-railed annulus.
\end{lemma}

\begin{proof}
Let $G$ be a planar graph of branch-width at least 11. Then by Theorem~\ref{thm:bnwdthcyl}, $G$ contains $(4,4)$-cylinder as a minor. By the definition of the minor relation, $G$ contains a $(4,4)$-railed annulus.
\end{proof}

\paragraph{Confluent paths} Let $G$ be a graph embedded in some surface $\Sigma$ 
and let $x\in V(G)$. We define a {\em disk around $x$}
as any open disk $\Delta_{x}$ with the property  that each point in $\Delta_{x}\cap G$ is  either  $x$ or belongs to the edges incident to $x$.
Let $P_{1}$ and $P_{2}$ be two edge-disjoint 
paths in $G$. 
We say that $P_{1}$ and $P_{2}$ are {\em confluent} if 
for every $x\in V(P_{1})\cap V(P_{2})$, that is not an endpoint 
of $P_{1}$ or $P_{2}$, and for every disk $\Delta_{x}$ around $x$, 
one of the connected components of 
the set $\Delta_{x}\setminus P_{1}$ does not contain any point of $P_{2}$.
We also say that a collection of paths is {\em confluent} if the paths in it are pairwise confluent. 

Moreover, given two edge-disjoint paths $P_{1}$ and $P_{2}$ in $G$ we say that a vertex $x\in V(P_{1})\cap V(P_{2})$ that is not an endpoint of $P_{1}$ or $P_{2}$ is an {\em overlapping vertex of $P_{1}$ and $P_{2}$} if there exists a $\Delta_{x}$ around $x$ such that both connected components of $\Delta_{x}\setminus P_{1}$ contain points of $P_{2}$. For a family of paths ${\cal P}$, a vertex $v$ of a path $P\in {\cal P}$ is called an {\em overlapping vertex of $P$} if there exists a path $P'\in {\cal P}$ such that $v$ is an overlapping vertex of $P$ and $P'$.

\section{Preliminary results on the confluency of paths}

\begin{lemma}\label{lem:pths}
Let $G$ be a graph and $v,v_{1},v_{2}\in V(G)$ such that there exist 
edge-disjoint paths $P_{1}$ and $P_{2}$ from $v$ to $v_{1}$ and $v_{2}$ respectively.
If the paths $P_{1}$ and $P_{2}$ are not well-arranged then there exist edge-disjoint paths $P_{1}'$ and $P_{2}'$ from $v$ to $v_{1}$ and $v_{2}$ respectively such that $E(P_{1}')\cup E(P_{2}')\subsetneq E(P_{1})\cup E(P_{2})$. 
\end{lemma}

\begin{proof}
Let $Z=V(P_{1})\cap V(P_{2})=\{v,u_{1},u_{2},\dots,u_{k}\}$, where $(v,u_{1},u_{2},\dots,u_{k})$ is the order that the vertices in $Z$ appear in $P_{1}$ and, $(v,u_{i_{1}},u_{i_{2}},\dots,u_{i_{k}})$ is the order that they appear in $P_{2}$.
As the paths are not well-arranged there exists $\lambda\in [k]$ such that $u_{\lambda}\neq u_{i_{\lambda}}$. Without loss of generality assume that $\lambda$ is the smallest such integer.
Without loss of generality 
assume also that $u_{\lambda}<u_{i_{\lambda}}$.
We define \begin{eqnarray*}
P_{1}' & = & P_{1}[v,u_{\lambda-1}]\cup P_{2}[u_{\lambda-1},u_{i_{\lambda}}]\cup P_{1}[u_{i_{\lambda}},v_{1}]\\ 
P_{2}' & = & P_{2}[v,u_{\lambda-1}]\cup P_{1}[u_{\lambda-1}, u_{\lambda}]\cup P_{2}[u_{\lambda},v_{2}].
\end{eqnarray*}
and observe that $P_{1}'$ and $P_{2}'$ satisfy the desired properties. (For an example, see Figure~\ref{fig:pths}).
\begin{figure}
\label{fig:pths}
\begin{center}
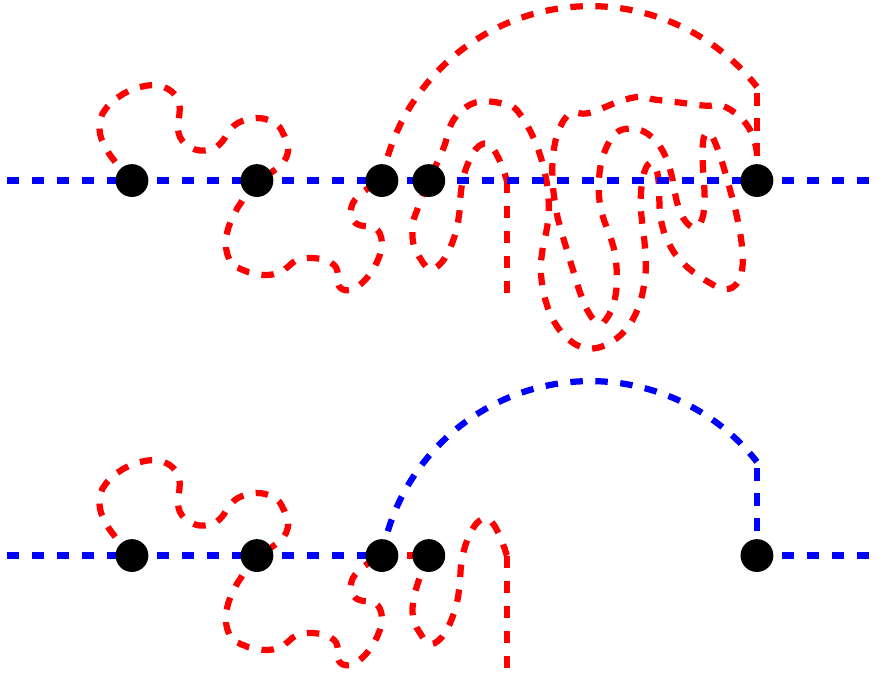
\caption{An example of the procedure in Lemma~\ref{lem:pths}}
\end{center}
\end{figure}
\end{proof}


Before proceeding to the statement and proof of the next proposition we need the following definition.
Given a collection of paths ${\cal P}$ in a graph $G$, we define the function $f_{\cal P}: \bigcup_{P\in {\cal P}}V(P) \rightarrow \mathbb{N}$ such that $f(x)$ is the number of pairs of paths $P,P'\in{\cal P}$ for which $x$ is an  overlapping vertex. Let 
$$\displaystyle g({\cal P})=\sum_{x\in \bigcup_{P\in {\cal P}}V(P)}f_{\cal P}(x).$$ 
Notice that $f(x)\geq 0$ for every $x\in \bigcup_{P\in {\cal P}}V(P)$ and thus $g({\cal P})\geq 0$. Observe also that $g({\cal P})=0$ if and only if ${\cal P}$ is a confluent collection of paths. 

Lemma~\ref{lem:pths} allows us to prove the main result of this section. We state the result for 
general surfaces as the proof for this more general setting does not have any essential difference than the case where $\Sigma$ is the  sphere $\Bbb{S}_{0}$.

\begin{proposition}
\label{tllds}
Let $r$ be a positive integer. If $G$ is a graph embedded in a surface $\Sigma$, $v,v_{1},v_{2},\dots,v_{r}\in V(G)$ and ${\cal P}$ 
is a collection of $r$ edge-disjoint paths from $v$ to $v_{1},v_{2},\dots, v_{r}$ in $G$, then $G$ contains a confluent collection ${\cal P}'$ of $r$ well-arranged edge-disjoint 
paths from $v$ to $v_{1},v_{2},\dots,v_{r}$ where $|{\cal P}'|=|{\cal P}|$
and such that $E(\bigcup_{P\in{\cal P}'}P)\subseteq E(\bigcup_{P\in{\cal P}}P)$.
\end{proposition}
\begin{proof}
Let $\hat{G}$ be the spanning subgraph of $G$ induced by the edges of the paths in ${\cal P}$ and let $G'$ be a minimal spanning subgraph of $\hat{G}$ that contains a collection of $r$ edge-disjoint paths from $v$ to $v_{1},v_{2},\dots, v_{r}$.
Let also ${\cal P}'$ be the collection of $r$ edge-disjoint paths from $v$ to $v_{1},v_{2},\dots, v_{r}$ in $G'$ for which $g({\cal P}')$ is minimum.
It is enough to prove that $g({\cal P}')=0$.

For a contradiction, we assume that $g({\cal P}')>0$ and we  prove that there exists a collection $\tilde{{\cal P}}$ of $r$ edge-disjoint paths from $v$ to $v_{1},v_{2},\dots, v_{r}$ in $G'$ such that $g(\tilde{P})<g({\cal P}')$. As $g({\cal P}')>0$, then there exists a path, say $P_{1}\in {\cal P}'$, that contains an overlapping vertex $u$. 
Let $z_{1}$ be the endpoint of $P_{1}$ which is different from $v$. 
Without loss of generality we may assume that $u$ is the overlapping vertex of $P$ that is closer to $z_{1}$ in $P$. 
Then there is a $(v,z_{2})$-path $P_{2}\in {\cal P}'$ such that $u$ is an overlapping vertex of $P_{1}$ and $P_{2}$.
Let $\tilde{P_{i}}=P_{3-i}[v,u]\cup P_{i}[u,z_{i}]$, $i\in[2]$ and $\tilde{P}=P$ for every $P\in {\cal P}'\setminus\{P_{1},P_{2}\}$.
As Lemma~\ref{lem:pths} and the edge-minimality of $G'$ imply that the paths $P_{1}$ and $P_{2}$ are well-arranged, we obtain that $\tilde{P}_{i}$ is a path from $v$ to $v_{i}$, $i\in [2]$. Let $\tilde{\cal P}$ be $\{\tilde{P}\mid P\in {\cal P}'\}$. It is easy to verify that $\tilde{\cal P}$ is a collection of $r$ edge-disjoint paths from $v$ to $v_{1},v_{2},\dots, v_{r}$.
We will now prove that $g(\tilde{{\cal P}})<g({\cal P}')$. 

First notice that if $x\neq u$, then $f_{\tilde{\cal P}}(x)=f_{{\cal P}'}(x)$. Thus, it is enough to prove that $f_{\tilde{\cal P}}(u)<f_{{\cal P}'}(u)$.
Observe that if $\{P,P'\}\subseteq {\cal P}'\setminus \{P_{1},P_{2}\}$ and $u$ is an overlapping vertex of $P$ and $P'$ then $u$ is also an overlapping vertex of $\tilde{P}$ and $P'$. Furthermore, while $u$ is an overlapping vertex in the case where $\{P,P'\}=\{P_{1},P_{2}\}$, it is not an overlapping vertex of $\tilde{P}_{1}$ and $\tilde{P}_{2}$. It remains to examine the case where $|\{P,P'\}\cap\{P_{1},P_{2}\}|=1$. In other words, we examine the case where one of the paths $P$ and $P'$, say $P'$, is $P_{1}$ or $P_{2}$, and $P\in {\cal P}'\setminus \{P_{1},P_{2}\}$. Let $\Delta_{u}$ be a disk around $u$ and $\Delta_{1},\Delta_{2}$ be the two distinct disks contained in the interior of $\Delta_{u}$ after removing $P$. 
We distinguish the following cases.\\

\noindent {\em Case 1.} $u$ is neither an overlapping vertex of $P_{1}$ and $P$, nor of $P_{2}$ and $P$ (see Figure~\ref{fig:pthcas1}). Then it is easy to see that the same holds for the pairs of paths $\tilde{P}_{1}$ and $P$ and, $\tilde{P}_{2}$ and $P$. Indeed, notice that for every $i\in [2]$, $P_{i}$ intersects exactly one of $\Delta_{1}$ and $\Delta_{2}$. Furthermore, as $u$ is an overlapping vertex of $P_{1}$ and $P_{2}$, both paths intersect the same disk. From the observation that $P_{1}\cup P_{2}=\tilde{P_{1}}\cup \tilde{P}_{2}$, we obtain that $u$ is neither an overlapping vertex of $\tilde{P}_{1}$ and $P$ nor of $\tilde{P}_{2}$ and $P$.\\

\begin{figure}\label{fig:pthcas1}
\begin{center}
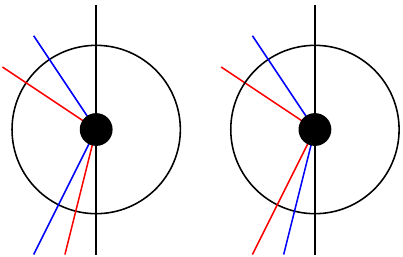
\end{center}
\caption{The paths $P$ (black), $P_{1}$ (red) and $P_{2}$ (blue) and the paths $\tilde{P}_{1}$ (blue) and $\tilde{P}_{2}$ (red).}
\end{figure}

\noindent {\em Case 2.} $u$ is an overlapping vertex of $P_{i}$ and $P$ but not of $P_{3-i}$ and $P$, $i\in [2]$ (see Figure~\ref{fig:pthcas2}).
Notice that exactly one of the following holds.
\begin{itemize}
\item $P_{i}[v,u]\cup P_{3-i}[v,u]$ intersects exactly one of the disks $\Delta_{1}$ or $\Delta_{2}$, say $\Delta_{1}$. Then $P_{i}[u,z_{i}]$ intersects $\Delta_{2}$ and $P_{3-i}[u,z_{3-i}]$ intersects $\Delta_{1}$. Therefore, it is easy to see that, $u$ is not an overlapping vertex of $P_{i}$ and $P$ anymore but becomes an overlapping vertex of $\tilde{P}_{3-i}$ and $P$.
\item $P_{i}[u,z_{i}]\cup P_{3-i}[u,z_{3-i}]$ intersects exactly one of the disks $\Delta_{1}$ or $\Delta_{2}$, say $\Delta_{1}$. Then $P_{i}[v,u]$ intersects $\Delta_{2}$ and $P_{3-i}[v,u]$ intersects $\Delta_{1}$. Therefore, it is easy to see that $u$ remains an overlapping vertex of $\tilde{P}_{i}$ and $P$ and does not become an overlapping vertex of $P_{3-i}$ and $P$.
\end{itemize}

\begin{figure}\label{fig:pthcas2}
\begin{center}
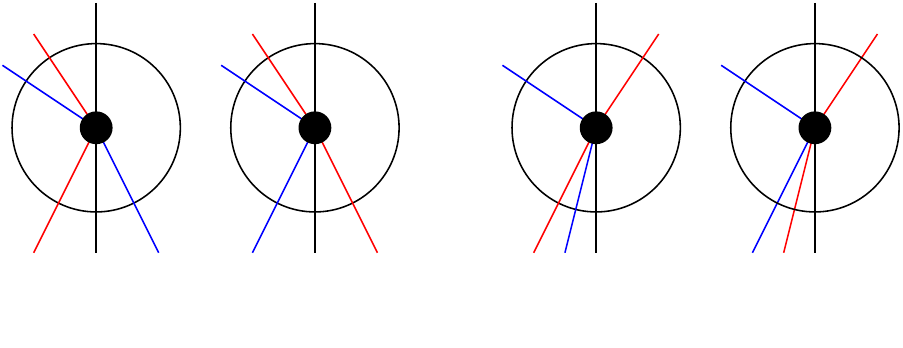
\end{center}
\caption{The paths $P$ (black), $P_{1}$ (red) and $P_{2}$ (blue) and the paths $\tilde{P}_{1}$ (blue) and $\tilde{P}_{2}$ (red).}
\end{figure}

\noindent {\em Case 3.} $u$ is an overlapping vertex of both $P_{1}$ and $P$ and, $P_{2}$ and $P$ (see Figure~\ref{fig:pthcas3}).
As above, exactly one of the following holds.
\begin{itemize}
\item $P_{1}[v,u]\cup P_{2}[v,u]$ intersects exactly one of the disks $\Delta_{1}$ or $\Delta_{2}$, say $\Delta_{1}$.
Then $P_{1}[u,z_{1}]\cup P_{2}[u,z_{2}]$ intersects $\Delta_{2}$. It follows that $u$ is an overlapping vertex of both $\tilde{P}_{1}$ and $P$ and, $\tilde{P}_{2}$ and $P$
\item $P_{1}[v,u]\cup P_{2}[u,z_{2}]$ intersects exactly one of the disks $\Delta_{1}$ or $\Delta_{2}$, say $\Delta_{1}$. Then $P_{1}[u,z_{1}]\cup P_{2}[v,u]$ intersects $\Delta_{2}$. It follows that $u$ is neither an overlapping vertex of $\tilde{P}_{1}$ and $P$ nor of $\tilde{P}_{2}$ and $P$.
\end{itemize}

\begin{figure}\label{fig:pthcas3}
\begin{center}
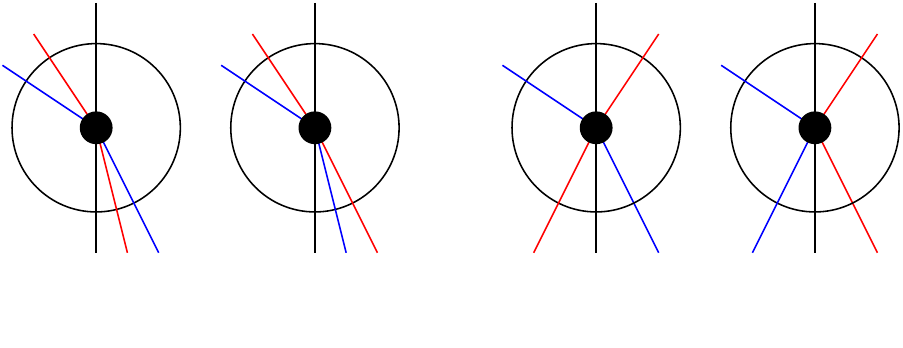
\end{center}
\caption{The paths $P$ (black), $P_{1}$ (red) and $P_{2}$ (blue) and the paths $\tilde{P}_{1}$ (blue) and $\tilde{P}_{2}$ (red).}
\end{figure}

From the above cases we obtain that $f_{\tilde{\cal P}}(u)<f_{{\cal P}'}(u)$ and therefore $g(\tilde{\cal P})<g({\cal P}')$, contradicting the choice of ${\cal P}'$. This completes the proof of the Proposition.
\end{proof}

\section{A decomposition theorem}

We prove the following decomposition theorem for $(K_{5},K_{3,3})$-immersion free graphs.

\begin{theorem}
\label{the:mainresult}
If $G$ is a graph not containing $K_{5}$ or $K_{3,3}$ as an immersion,
then $G$ can be constructed by applying consecutive $i$-edge sums, for $i\leq 3$, to graphs that either are sub-cubic or have branch-width at most 10.
\end{theorem}

\begin{proof}
Observe first that a $(K_{5},K_{3,3})$-immersion-free graph
is also $(K_{5},K_{3,3})$-topological-minor-free, therefore, from Kuratowski's
theorem, $G$ is planar. 
Applying Lemma~\ref{frloo}, we may assume that $G$ is a  
$(K_{5},K_{3,3})$-immersion-free graph $G$ without 
any internal $i$-edge cut, $i\in [3]$. It is now enough to prove that $G$ is either planar sub-cubic
or has branch-width at most 10. For a contradiction, we 
assume that $\bw(G)\geq 11$ and that 
$G$ contains some vertex $v$ of degree $\geq 4$. Our aim is to prove that $G$ 
contains $K_{3,3}$ as an immersion. First, let $G^{s}$ be the graph obtained from $G$ after subdividing all of its edges once. Notice that 
$G^{s}$ contains $K_{3,3}$ as an immersion if and only if $G$ contains $K_{3,3}$ as an immersion. Hence, from now on, we want to find $K_{3,3}$ in $G^{s}$ as an immersion.

From Lemma~\ref{lem:brnwdthranl} $G$, and thus $G^{s}$, contains a $(4,4)$-railed annulus as a subgraph.
Observe then that $G^{s}$ also contains as a subgraph a $(2,4)$-railed annulus such that the vertex $v$ of 
degree $\geq 4$ does not belong in the annulus between its cycles (Figure~\ref{edgsmbef} depicts 
the case where $v$ is inside the annulus between the second and the third cycle). We denote by $C_{1}$ and $C_{2}$ 
the nested cycles and by $R_{1},R_{2},R_{3}$ and $R_{4}$ the rails of the above $(2,4)$-railed annulus. Let $A$ 
be the annulus  between $C_{1}$ and $C_{2}$. Without loss of generality we may assume that $C_{1}$ 
separates $v$ from $C_{2}$ and that $A$ is edge-minimal, that is, there is no other annulus $A'$ such that $|E(A')|<|E(A)|$ and $A'\subseteq A$.
\begin{figure}[h]
\begin{center}
\scalebox{.42}{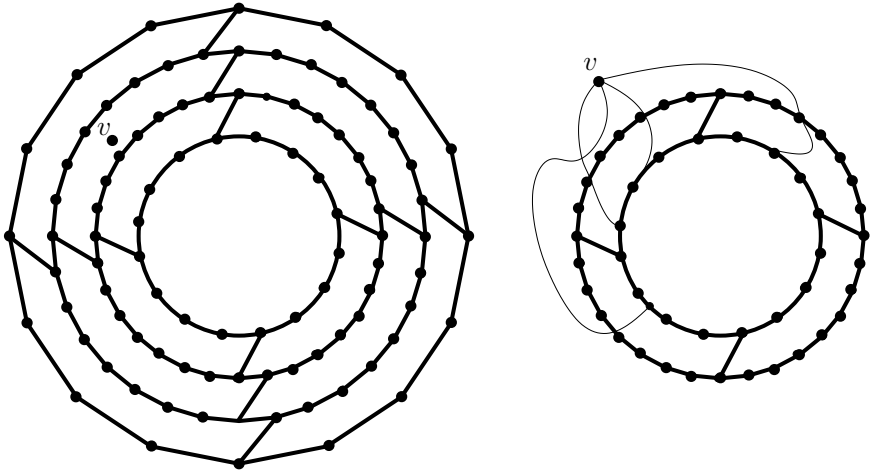} 
\end{center}
\caption{The $(4,4)$-railed annulus and the vertex $v$}
\label{edgsmbef}
\end{figure}

Let now $G_{1},G_{2},\dots,G_{p}$ be the connected components of $A\setminus (C_{1}\cup C_{2})$.

\begin{claim}\label{clm1}
For every $i\in[p]$ and every $j\in [2]$, $|N_{G^{s}}(V(G_{i}))\cap V(C_{j})|\leq 1$.
\end{claim}

\begin{proof}[Proof of Claim~\ref{clm1}]
Indeed, assume the contrary.
Then there is a cycle $C_{j}'$ such that $C_{j}'$ and $C_{j\mod 2 +1}$ define an annulus $A'$ with $A'\subseteq A$ and $|E(A')|< |E(A)|$, a contradiction to the edge-minimality of the annulus $A$.
\end{proof}

For every $l\in [p]$, we denote by $u^{l}_{1}$ and $u^{l}_{2}$ the unique neighbor of $G_{k}$ in $C_{1}$ and $C_{2}$ respectively (whenever they exist). We call the connected components that have both a neighbor in $C_{1}$ and a neighbor in $C_{2}$ {\em substantial}. 
Let ${\cal C}=\{\widehat{G}_{i}=G[V(G_{i})\cup \{u^{i}_{1},u^{i}_{2}\}]\mid G_{i}\text{ is a substantial connected component}\}$. That is, ${\cal C}$ is the set of graphs induced by the substantial connected components and their neighbors in the cycles $C_{1}$ and $C_{2}$. Note that every edge of $G$ has been subdivided in $G^{s}$ and thus every edge $e\in G$ for which $e\cap C_{1}\neq \emptyset$ and $e\cap C_{2}\neq \emptyset$ corresponds to a substantial connected component in ${\cal C}$.

We now claim that there exist four confluent edge-disjoint paths $P_{1},P_{2},P_{3}$ and $P_{4}$ from $v$ to $C_{2}$ in $G^{s}$.
This follows from the facts that $G^{s}$ does not contain an internal $i$-edge cut, $C_{2}$ contains at least $4$ vertices, $\deg_{G^{s}}(v)\geq 4$ combined with Observation~\ref{obs:mngr}. Moreover, from Proposition~\ref{tllds}, we may assume that $P_{1},P_{2},P_{3}$ and $P_{4}$ are confluent. 


Let $P_{i}'$ be the subpath $P_{i}[v,v_{i}]$ of $P_{i}$, where $v_{i}$ is the vertex in $V(P_{i})\cap V(C_{2})$ whose distance from $v$ in $P_{i}$ is minimum, $i\in [4]$. Recall that all edges of $G$ have been subdivided in $G^{s}$. This implies that there exist four (possibly not disjoint) graphs in ${\cal C}$, say $\widehat{G}_{1},\widehat{G}_{2},\widehat{G}_{3}$ and $\widehat{G}_{4}$ such that $v_{i}=u^{i}_{2}$, $i\in [4]$. 
We distinguish two cases.\\

\noindent{\em Case 1.} The graphs $\widehat{G}_{1},\widehat{G}_{2},\widehat{G}_{3}$ and $\widehat{G}_{4}$ are vertex-disjoint.\\
This implies that the endpoints of $P_{1}',P_{2}',P_{3}'$ and $P_{4}'$ are disjoint. Let $G'$ be the graph induced by the cycles $C_{1},C_{2}$ and the paths $P_{1}',P_{2}',P_{3}',P_{4}'$
and let $\widehat{P}_{1},\widehat{P}_{2},\widehat{P}_{3}$ and $\widehat{P}_{4}$ be confluent edge-disjoint paths from $v$ to $u_{2}^{1}$, $u_{2}^{2}$, $u_{2}^{3}$ and $u_{2}^{4}$ in $G'$ such that
\begin{enumerate}[(i)]
\item \label{prop1} $\sum\{e\mid e\in \bigcup_{i\in[4]}E(\widehat{P}_{i})\setminus E(A)\}$ is minimum, that is, the number of the edges of the paths that is outside of $A$ is minimum, and
\item \label{prop2} subject to~\ref{prop1}, $\sum\{e\mid e\in \bigcup_{i\in[4]}E(\widehat{P}_{i})\}$ is minimum.
\end{enumerate}
Let also $\widehat{G}$ be the graph induced by $C_{1}$, $C_{2}$, $\widehat{P}_{1}$, $\widehat{P}_{2}$, $\widehat{P}_{3}$ and $\widehat{P}_{4}$. From now on we work towards showing that $\widehat{G}$ contains $K_{3,3}$ as an immersion. For every $i\in [4]$ we call a connected component of $\widehat{P}_{i}\cap C_{1}$ non-trivial if it contains at least an edge.

\begin{claim}\label{clm2}
For every $i\in [4]$, $\widehat{P}_{i}\cap C_{1}$ contains at most one non-trivial connected component $Q_{i}$ and $u_{1}^{i}$ is an endpoint of $Q_{i}$.
\end{claim}

\begin{proof}[Proof of Claim~\ref{clm2}]
First, notice that any path from $v$ to $v_{i}$ in $\widehat{G}$ contains $u_{1}^{i}$ and thus, $u_{1}^{i}\in V(\widehat{P}_{i})$. 
Observe now that $\widehat{P}_{i}[u_{1}^{i},u_{2}^{i}]$ is a subpath of $\widehat{P}_{i}$ whose internal vertices do not belong to $C_{1}$, thus if $u_{1}^{i}$ belongs to a non-trivial connected component $Q_{i}$ of $\widehat{P}_{i}\cap C_{1}$, then $u_{1}^{i}$ is an endpoint of $Q_{i}$.
We will now prove that any non-trivial connected component of $\widehat{P}_{i}\cap C_{1}$ contains $u_{1}^{i}$.
Assume in contrary that there exists a non-trivial connected component $P$ of $\widehat{P}_{i}\cap C_{1}$ that does not contain $u_{1}^{i}$. Let $u$ be the endpoint of $P$ for which $\dist_{\widehat{P}_{i}}(u,u_{1}^{i})$ is minimum. 
Let also $u'$ be the vertex in $\widehat{P}_{i}[u,u_{1}^{i}]\cap C_{1}$ such that $\dist_{\widehat{P}_{i}}(u,u')$ is minimum. 
Let $P'$ be the subpath of $C_{1}$ with endpoints $u,u'$ such that $\widehat{P}_{i}[u,u']\cup P'$ is a cycle $C$ with $C\cap P=\{u\}$. We further assume that the interior of $\widehat{P}_{i}[u,u']\cup P'$ is the open disk that does not contain any vertices of $\widehat{P}_{i}$.  We will prove that for every path $\widehat{P}_{j}$, $j\in[4]$ $\widehat{P}_{j}\cap P' \subseteq \{u,u'\}$. As this trivially holds for $j=i$ we will assume that $j\neq i$.
 Observe that, for every $j\in [4]$, $\widehat{P}_{j}[v,u_{1}^{j}]\cap A\subseteq C_{1}$ as for every connected component $H$ of $A\setminus (C_{1}\cup C_{2})$ it holds that $|N_{G^{s}}(V(H))\cap V(C_{j})|\leq 1$. Furthermore, observe that $\widehat{P}_{i}[u,u']\cup P'$ is a separator in $\widehat{G}$. 
This implies that $v$ does not belong to the interior of $\widehat{P}_{i}[u,u']\cup P'$. Thus, if there is a vertex $z$ such that $z\in \widehat{P}_{j}\cap (P'\setminus \{u,u'\})$, $j\neq i$ there is a vertex $z'\in \widehat{P}_{j}\cap \widehat{P}_{i}[u,u']$, a contradiction to the confluency of the paths. We may then replace $\widehat{P}_{i}[u,u']$ by $P'$, a contradiction to~\ref{prop1}. 
\end{proof}

We denote by $v_{i}$ the endpoint of $Q_{i}$ that is different from $u_{1}^{i}$ if $Q_{i}$ is a non-trivial connected component of $\widehat{P}_{i}\cap C_{1}$, $i\in [4]$. Observe that $\widehat{P}_{i}=\widehat{P}_{i}[v,v_{i}]\cup Q_{i} \cup \widehat{P}_{i}[u_{1}^{i},u_{2}^{i}]$, where we let $Q_{i}=\emptyset$ in the case where $\widehat{P}_{i}\cap C_{1}$ is edgeless, $i\in [4]$.
We denote by $T_{i}$ the subpath of $C_{1}$ with endpoints $u_{1}^{i}$ and $u_{1}^{i\mod4+1}$ such that $T_{i}\cap \{\{u_{1}^{1},u_{1}^{2},u_{1}^{3},u_{1}^{4}\}\setminus \{u_{1}^{i},u_{1}^{i\mod4+1}\}\}=\emptyset$, $i\in [4]$.
From the confluency of the paths $\widehat{P}_{i}$ and the fact that $u_{1}^{i}$ is an endpoint of $Q_{i}$ it follows that $Q_{i}\subseteq T_{i}$ or $Q_{i}\subseteq T_{i-1}$, $i\in[4]$ where $T_{i-1}=T_{3+i\mod 4}$ if $i-1\notin [4]$. 

\begin{claim}\label{clm3}
There exists an $i_{0}\in [4]$ such that $T_{i_{0}}\cap (Q_{i_{0}},Q_{i_{0}\mod4+1})\neq T_{i_{0}}$. 
\end{claim}

\begin{proof}[Proof of Claim~\ref{clm3}]
Towards a contradiction assume that for every $i\in [4]$, it holds that $T_{i}\cap (Q_{i},Q_{i\mod4+1})= T_{i}$.
It follows that either $Q_{i}=T_{i}=\widehat{P}_{i}[v_{i},v_{1}^{i}]$, $i\in [4]$ or $Q_{i\mod4+1}=T_{i}$, $i\in [4]$. 
Notice then that either $v_{i}=u_{1}^{i\mod 4+1}$, $i\in [4]$ or $v_{i\mod4+1}=u_{1}^{i}$, $i\in [4]$ respectively.  
Then, we let $\tilde{P}_{i\mod4+1}=\widehat{P}_{i}[v,v_{i}]\cup\widehat{P}_{i\mod4+1}[u_{1}^{i\mod4+1},u_{2}^{i\mod4+1}]$ or $\tilde{P}_{i}=\widehat{P}_{i\mod4+1}[v,v_{i\mod4+1}]\cup \widehat{P}_{i}[u_{1}^{i},u_{2}^{i}]$, $i\in [4]$ respectively.
Notice that the paths $\tilde{P}_{1},\tilde{P}_{2},\tilde{P}_{3}$ and $\tilde{P}_{4}$ are confluent edge-disjoint paths from $v$ to $u_{2}^{1},u_{2}^{2},u_{2}^{3}$ and $u_{2}^{4}$ such that $\cup_{i\in[4]}\tilde{P}_{i}$ is a proper subgraph of $\cup_{i\in[4]}\widehat{P}_{i}$. Therefore, we have that $\sum\{e\mid e\in \bigcup_{i\in[4]}E(\tilde{P}_{i})\}<\sum\{e\mid e\in \bigcup_{i\in[4]}E(\widehat{P}_{i})\}$, a contradiction to~\ref{prop2}. 
\end{proof}

It is now easy to see that $\widehat{G}$, and thus $G$, contains $K_{3,3}$ as an immersion. 
Indeed, first remove all edges of $C_{1}\setminus T_{i_{0}}$ that do not belong to any path $\widehat{P}_{i}$, $i\in [4]$. Then lift the paths $\widehat{P}_{i}$ to a single edge where $i\neq i_{0},i_{0}\mod 4+1$.
Now let $u_{i_{0}}$ ($u_{i_{0}\mod4+1}$ respectively) be the vertex of $T_{i_{0}}$ that belongs to $\widehat{P}_{i_{0}}$ ($\widehat{P}_{i_{0}\mod4+1}$ respectively) whose distance from $v$ in $\widehat{P}_{i_{0}}$ ($\widehat{P}_{i_{0}\mod4+1}$ respectively) is minimum and lift the paths $\widehat{P}_{i_{0}}[v,u_{i_{0}}]$ and $\widehat{P}_{i_{0}\mod4+1}[v,u_{i_{0}\mod4+1}]$ to single edges.
Notice now that $\widehat{G}$ contains the graph $H_{2}$ depicted in Figure~\ref{3topmin} as an immersion. Thus, we get that $\widehat{G}$ contains $K_{3,3}$ as an immersion.\\

\begin{figure}[h]
\begin{center}
\scalebox{.79}{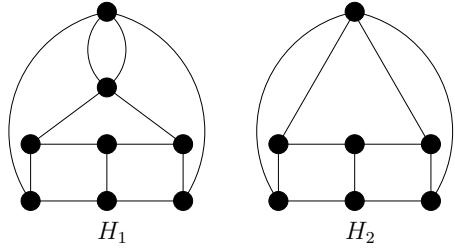} 
\end{center}
\caption{The graphs $H_{1}$ and $H_{2}$}
\label{3topmin}
\end{figure}

\noindent{\em Case 2.} There exist $i_{1},i_{2}\in [4]$ such that $\widehat{G}_{i_{1}}$ and $\widehat{G}_{i_{2}}$ are not vertex-disjoint.\\
Let $G^{\mu}$ be the graph induced by the cycles $C_{1}$ and $C_{2}$ and the graphs in ${\cal C'}$. We will show that $G^{\mu}$ contains $K_{3,3}$ as an immersion. First recall that the common vertices of $\widehat{G}_{i_{1}}$ and $\widehat{G}_{i_{2}}$ lie in at least one of the cycles $C_{1}$ and $C_{2}$. Without loss of generality assume that they have a common vertex in $C_{1}$. Recall that, as every edge of $G$ has been subdivided in $G^{s}$, there does not exist an edge $e\in G^{s}$ such that $e\cap C_{j}\neq\emptyset$, $j\in[2]$. This observation and the fact that there exist four rails between $C_{1}$ and $C_{2}$ imply that there exist at least four graphs in ${\cal C}'$ that are vertex-disjoint. It follows that there exist three vertex-disjoint graphs, say $\widehat{G}_{i_{3}},\widehat{G}_{i_{4}},\widehat{G}_{i_{5}}$, in ${\cal C}'$ with the additional properties that 
$\widehat{G}_{i_{2+r}}\cap \widehat{G}_{i_{1}} \cap C_{1}=\emptyset$, $r\in [3]$ and that at most one of the $\widehat{G}_{i_{3}}, \widehat{G}_{i_{4}},\widehat{G}_{i_{5}}$ has a common vertex with one of the $\widehat{G}_{i_{1}},\widehat{G}_{i_{2}}$. Note here that none of the $\widehat{G}_{i_{3}}, \widehat{G}_{i_{4}},\widehat{G}_{i_{5}}$ can have a common vertex with one of the $\widehat{G}_{i_{1}},\widehat{G}_{i_{2}}$ in $C_{2}$, in the case where $\widehat{G}_{i_{1}}\cap \widehat{G}_{i_{2}}\cap C_{2}\neq\emptyset$. It is now easy to see that $G^{\mu}$ contains $H_{1}$ or ($H_{2}$ respectively) depicted in Figure~\ref{3topmin} as a topological minor when $\widehat{G}_{i_{1}}\cap \widehat{G}_{i_{2}}\cap C_{2}\neq\emptyset$ ($\widehat{G}_{i_{1}}\cap \widehat{G}_{i_{2}}\cap C_{2}=\emptyset$ respectively).
Observe now that $H_{1}$ contains $H_{2}$ as an immersion. Moreover, notice that $H_{2}$ contains $K_{3,3}$ as an immersion. Thus $G^{\mu}$, and therefore $G^{s}$ and $G$, contain $K_{3,3}$ as an immersion, a contradiction.
\end{proof}

\begin{remark}
It is easy to verify that our results hold for both  the weak and strong immersion relations.
\end{remark}

\begin{figure}[h]
\begin{center}
\scalebox{.2}{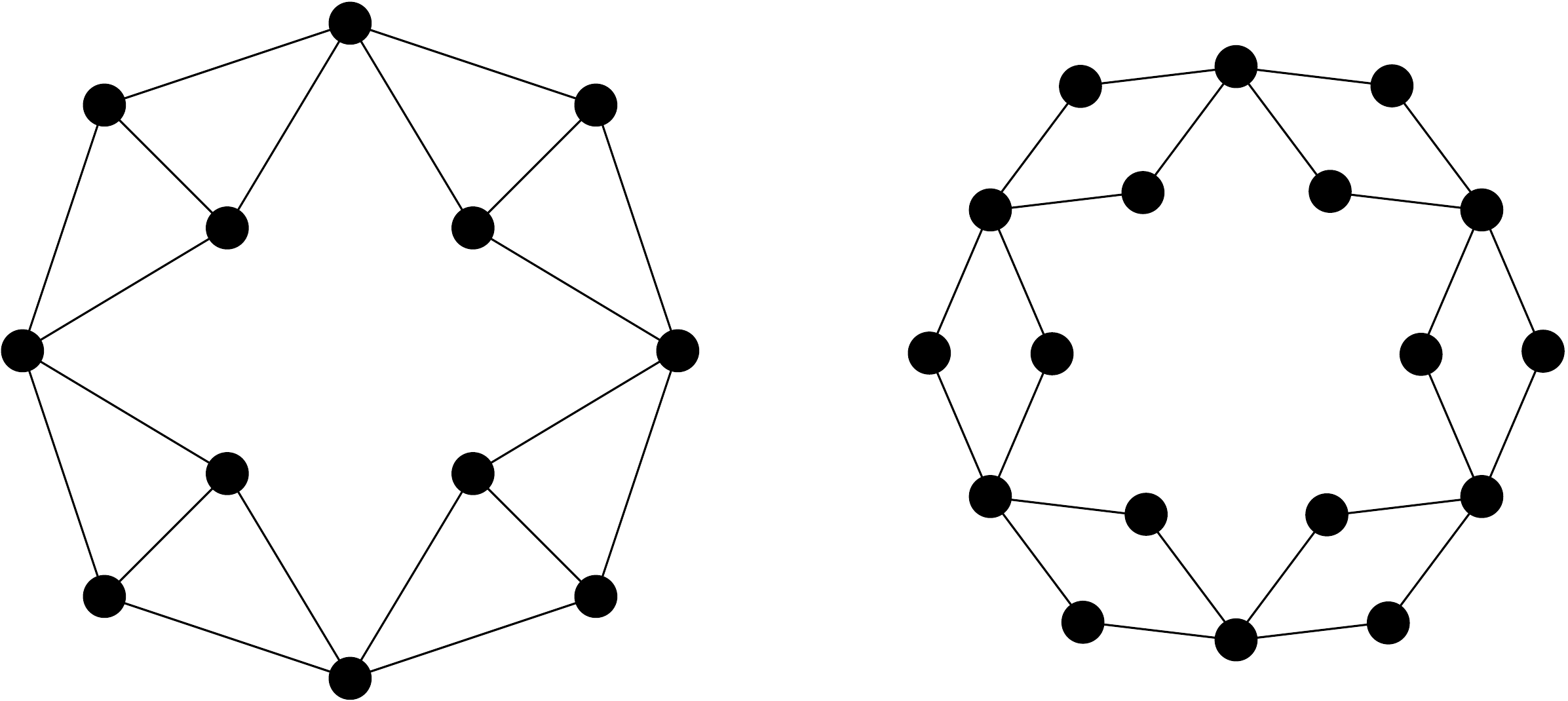} 
\end{center}
\caption{Simple non-sub-cubic graphs of branch-width 3 without $K_{5}$ or $K_{3,3}$ as immersions.}
\label{3topmin}
\end{figure}

We believe that the upper bound on the branch-width 
of the building blocks of Theorem~\ref{the:mainresult}
can be further reduced, especially if we restrict ourselves
to simple graphs. There are infinite such graphs that are not sub-cubic
and have  branch-width 3; some of them are depicted in Figure~\ref{3topmin}.
However, we have not been able to 
find any simple non-sub-cubic graph of branch-width greater than 3 that does not contain
$K_{5}$ or $K_{3,3}$ as an immersion.

\bibliography{complete}
\bibliographystyle{plain}
\end{document}